\documentclass[12pt]{article}
\usepackage{amsfonts}
= 0.0in \topmargin = 0.3in \oddsidemargin=0.1in
\def\Dj{\hbox{D\kern-.73em\raise.30ex\hbox{-}
\raise-.30ex\hbox{}}}
\def\dj{\hbox{d\kern-.33em\raise.80ex\hbox{-}
\raise-.80ex\hbox{\kern-.40em}}}
\usepackage{epsfig}
\usepackage{amsmath,amsthm,amsfonts,amssymb,amscd}
\usepackage{graphicx,ifpdf}

\begin{document}
\baselineskip=0.30in

\newtheorem{lem}{Lemma}[section]
\newtheorem{thm}[lem]{Theorem}
\newtheorem{cor}[lem]{Corollary}
\newtheorem*{prop}{Proposition}
\newtheorem{con}[lem]{Conjecture}
\newtheorem{rem}[lem]{Remark}
\newtheorem{defi}[lem]{Definition}
\renewcommand\baselinestretch{1.2}
\def\pf{\noindent {\it Proof.} }
\def\qed{\hfill \rule{4pt}{7pt}}


\renewcommand\thefootnote{\fnsymbol{footnote}}

\begin{center} {\Large \bf On the minimal energy of conjugated unicyclic graphs with maximum degree at most
3}
 \end{center}

\baselineskip=0.20in

\vspace{2mm}

\baselineskip=0.20in
\begin{center}
 { \small   Hongping Ma$^{1}$, Yongqiang Bai$^{1}$\footnote{Corresponding author.}, Shengjin Ji$^{2}$ \\[5pt]
\small $^{1}$ School of Mathematics and Statistics, Jiangsu Normal University,\\ Xuzhou 221116, China\\
 \small $^{2}$ School of Science, Shangdong University of Technology, \\ Zibo 255049, China\\
 \small Email:  hpma@163.com, bmbai@163.com, jishengjin2013@163.com }
\end{center}


\begin{abstract}
 The energy of a graph $G$, denoted by $E(G)$, is defined as the sum of the absolute values of all eigenvalues of $G$.
 Let $n$ be an even number and $\mathbb{U}_{n}$ be the set of all conjugated unicyclic graphs of order $n$ with maximum degree at most $3$.
 Let $S_n^{\frac{n}{2}}$ be the radialene graph obtained by attaching a pendant
 edge to each vertex of the cycle $C_{\frac{n}{2}}$.
   In [Y. Cao et al., On the minimal energy of unicyclic H\"{u}ckel molecular graphs possessing Kekul\'{e} structures,
 Discrete Appl. Math. 157 (5) (2009), 913--919], Cao et al. showed
 that if $n\geq 8$, $S_n^{\frac{n}{2}}\ncong G\in \mathbb{U}_{n}$ and the girth of $G$ is not divisible by $4$, then $E(G)>E(S_n^{\frac{n}{2}})$.
 Let $A_n$ be the unicyclic graph obtained by attaching a $4$-cycle to one of the two leaf  vertices of the path $P_{\frac{n}{2}-1}$ and a pendent edge to
 each other vertices of $P_{\frac{n}{2}-1}$.
 In this paper, we prove that $A_n$ is the unique unicyclic graph in
 $\mathbb{U}_{n}$ with minimal energy.
\\[2mm]
\end{abstract}

\baselineskip=0.27in

\section{Introduction}

 Let $G$ be a simple graph with $n$ vertices and $A(G)$ the adjacency matrix of
 $G$. The eigenvalues $\lambda_{1}, \lambda_{2},\ldots, \lambda_{n}$ of $A(G)$  are said to be the eigenvalues of the graph $G$.
The energy of $G$ is defined as
$$E=E(G)=\sum_{i=1}^{n}|\lambda_{i}|.$$
This concept was intensively studied in chemistry, since it can be
used to approximate the total $\pi$-electron energy of a molecular.
Further details on the mathematical properties and chemical
 applications of $E(G)$, see the recent book \cite{LSG}, reviews \cite{G2,GLZ}, and papers \cite{BB,DM,DMG,GFAHG,LSWL,MMZ,Z}.

\begin{figure}[ht]
\centering
  \setlength{\unitlength}{0.05 mm}%
  \begin{picture}(1784.8, 562.2)(0,0)
  \put(0,0){\includegraphics{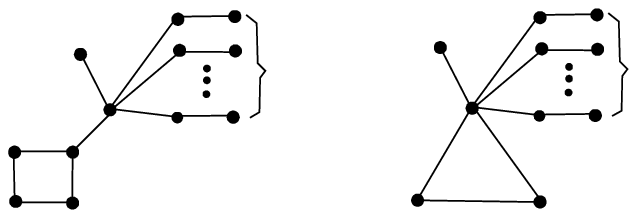}}
  \put(73.92,38.20){\fontsize{8.53}{10.24}\selectfont \makebox(150.0, 60.0)[l]{$U_1$\strut}}
  \put(1294.83,382.88){\fontsize{8.53}{10.24}\selectfont \makebox(450.0, 60.0)[l]{$\frac{n}{2}-2$\strut}}
  \put(561.90,383.72){\fontsize{8.53}{10.24}\selectfont \makebox(450.0, 60.0)[l]{$\frac{n}{2}-3$\strut}}
  \put(944.44,58.39){\fontsize{8.53}{10.24}\selectfont \makebox(150.0, 60.0)[l]{$U_2$\strut}}
  \end{picture}%
\caption{The conjugated unicyclic graphs with minimal and
second-minimal energy.}\label{fig-minimal}
\end{figure}

 One of the fundamental question that is encountered in the study of graph energy is which graphs (from a given class) have minimal and maximal energies.
 A large of number of papers were published on such extremal problems, especially for various subclasses of trees and unicyclic graphs, see Chapter 7 in \cite{LSG}.
 A conjugated unicyclic graph is a connected graph with one unique cycle that has a perfect matching.
 The problem of determining the conjugated unicyclic graph with minimal energy has been considered in \cite{LZZ,WCL}, and
 Li et al. \cite{LZZ} proved that the conjugated unicyclic graph of order
(even) $n$ with minimal energy  is $U_1$ or $U_2$, as shown in
Figure \ref{fig-minimal}.
 It has been shown that $E(U_1)<E(U_2)$ by Li and Li \cite{LL}.
 Recently, results on ordering of conjugated unicyclic graphs by minimal energies have been extended in \cite{W,
 Z}. In particular, $U_2$ is unique conjugated unicyclic graph of order $n$ with second-minimal
 energy.

 The degree of a vertex $v$ in a graph $G$ is denoted by $d_G(v)$. Denote by
$\Delta$ the maximum degree of a graph. From now on, let $n$ be an
even number. Let $\mathbb{U}_{n}$ be the set of all conjugated
unicyclic graphs of order $n$ with $\Delta\leq 3$.
 Let $G\in \mathbb{U}_{n}$, the length of the unique cycle of $G$ is denoted by $g(G)$, or simply
 $g$, and the unique cycle of $G$ is denoted by $C_g(G)$, or simply $C_g$.
 Let $S_n^{\frac{n}{2}}$ be the radialene graph obtained by attaching a pendant
 edge to each vertex of the cycle $C_{\frac{n}{2}}$. Wang et al. \cite{WCZL}
 showed the following results: Assume that $n\geq 6$ and $S_n^{\frac{n}{2}}\ncong G\in
 \mathbb{U}_{n}$. Then if one of the following conditions holds: (i)
 $\frac{n}{2}\equiv g\equiv 1$ (mod $2$) and $g\leq\frac{n}{2}$, (ii) $g\not\equiv \frac{n}{2}\equiv 0$ (mod
 $4$), (iii) $\frac{n}{2}\equiv g\equiv  2$ (mod $4$), and
 $g\leq\frac{n}{2}$, then $E(G)>E(S_n^{\frac{n}{2}})$. Y. Cao et al. \cite{CLLZ} improved the
 above results by proving the following Lemma.

\begin{lem} \textnormal{\cite{CLLZ}}\label{lem g-non-4-multiple}
 If $n\geq 8$, $S_n^{\frac{n}{2}}\ncong G\in \mathbb{U}_{n}$ with $g\not\equiv  0$ \textnormal{(mod
 $4$)}, then $E(G)>E(S_n^{\frac{n}{2}})$.
\end{lem}

\begin{figure}[ht]
\centering
   \setlength{\unitlength}{0.05 mm}%
  \begin{picture}(2792.5, 525.4)(0,0)
  \put(0,0){\includegraphics{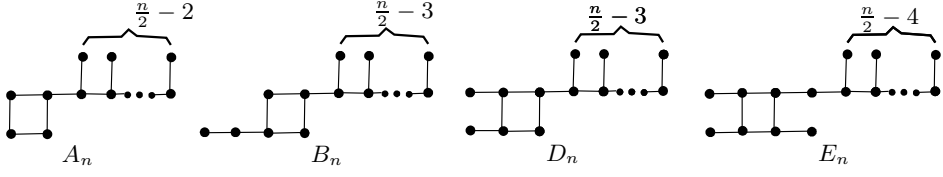}}
  \put(185.67,39.45){\fontsize{8.53}{10.24}\selectfont \makebox(150.0, 60.0)[l]{$A_n$\strut}}
  \put(374.88,419.32){\fontsize{8.53}{10.24}\selectfont \makebox(450.0, 60.0)[l]{$\frac{n}{2}-2$\strut}}
  \put(847.63,38.20){\fontsize{8.53}{10.24}\selectfont \makebox(150.0, 60.0)[l]{$B_n$\strut}}
  \put(1012.53,423.62){\fontsize{8.53}{10.24}\selectfont \makebox(450.0, 60.0)[l]{$\frac{n}{2}-3$\strut}}
  \put(1472.33,42.99){\fontsize{8.53}{10.24}\selectfont \makebox(150.0, 60.0)[l]{$D_n$\strut}}
  \put(1578.87,409.06){\fontsize{8.53}{10.24}\selectfont \makebox(450.0, 60.0)[l]{$\frac{n}{2}-3$\strut}}
  \put(1578.87,409.06){\fontsize{8.53}{10.24}\selectfont \makebox(450.0, 60.0)[l]{$\frac{n}{2}-3$\strut}}
  \put(2196.00,41.96){\fontsize{8.53}{10.24}\selectfont \makebox(150.0, 60.0)[l]{$E_n$\strut}}
  \put(2302.54,408.04){\fontsize{8.53}{10.24}\selectfont \makebox(450.0, 60.0)[l]{$\frac{n}{2}-4$\strut}}
  \end{picture}%
 \caption{Four graphs in $\mathbb{U}_{n}$.}\label{fig-minimal-degree-3}
\end{figure}

 Let $A_n$, $B_n$, $D_n$ and $E_n$ be the graphs shown in Figure \ref{fig-minimal-degree-3}.

 In this paper, we completely characterize the graph with minimal
 energy in $\mathbb{U}_{n}$ by showing the following result.


\begin{thm}\label{main thm}
 $A_n$ is the unique unicyclic graph in $\mathbb{U}_{n}$ with minimal
 energy for $n\geq 8$.
\end{thm}

\section{Preliminaries}

 In this section, we first introduce some notations and properties which are need in the
 sequel. Then we give some results on the energies of graphs $A_n$, $B_n$,
 $D_n$, $E_n$ and $S_n^{\frac{n}{2}}$.

 Let $G$ be a graph of order $k$. The characteristic polynomial of $A(G)$ is also called the characteristic polynomial of $G$, denoted by
 $\phi(G,x)=\mbox{det}(xI-A(G))=\sum_{i=0}^{k}a_i(G)x^{k-i}$. Using
 these coefficients of $\phi(G,x)$, the energy of $G$ can be expressed as the Coulson integral formula
 \cite{GP}:
  \begin{eqnarray}
 E(G)=\frac{1}{2\pi}{\large\int}_{-\infty}^{+\infty}\frac{1}{x^{2}}\ln\left[\left(\sum\limits_{i=0}^{\lfloor
\frac{k}{2}\rfloor}
(-1)^ia_{2i}(G)x^{2i}\right)^2+\left(\sum\limits_{i=0}^{\lfloor
\frac{k}{2}\rfloor} (-1)^ia_{2i+1}(G)x^{2i+1}\right)^2\right]dx.
\label{energy-1}
\end{eqnarray}
 Write $b_i(G)=|a_i(G)|$. Clearly, $b_0(G)=1$, $b_1(G)=0$, and $b_2(G)$ equals the number of edges of $G$. For unicyclic graphs or
 bipartite graphs, it can be shown \cite{GP,H} that
\begin{eqnarray}
 E(G)=\frac{1}{2\pi}{\large\int}_{-\infty}^{+\infty}\frac{1}{x^{2}}\ln\left[\left(\sum\limits_{i=0}^{\lfloor
\frac{k}{2}\rfloor}
b_{2i}(G)x^{2i}\right)^2+\left(\sum\limits_{i=0}^{\lfloor
\frac{k}{2}\rfloor} b_{2i+1}(G)x^{2i+1}\right)^2\right]dx.
\label{energy-2}
\end{eqnarray}
By formula \eqref{energy-2}, it is convenient to introduce the
following quasi-order relation \cite{SSGL}: if $G_1$ and $G_2$ are
two unicyclic or bipartite graphs with $k$ vertices, then
\[ G_1\succeq G_2 \Leftrightarrow b_i(G_1)\geq b_i(G_2) \mbox{ for
all } i=2, \ldots, k.
\]
 If  $G_1\succeq G_2$ and there exists some $j$ such that $b_j(G_1) > b_j(G_2)$, then we write $G_1\succ G_2$.
 Clearly,  $G_1\succeq G_2\Rightarrow  E(G_1)\geq E(G_2)$, and $G_1\succ G_2\Rightarrow  E(G_1)>E(G_2)$.
 It is known \cite{CDS} that  $b_{2i+1}(G)=0$ for a bipartite graph $G$,
 and $b_{2i}(G)$ equals the number of $i$-matchings of $G$ if $G$ is a tree.

\begin{lem} \textnormal{\cite{LZ}}\label{lem delete-edge}
 Let $G$ be a graph whose components are all trees except at most one being a unicyclic graph.

 (1) If $G$ contains exactly one cycle $C_g$, and $uv$ is an edge on this cycle, then
\[
\begin{array}{lll}
b_i(G)=b_i(G-uv) + b_{i-2}(G-u-v)-2b_{i-g}(G-C_g)\ \mbox{if}\
g\equiv 0
 \  (\textnormal{mod} \ 4),\\
b_i(G)=b_i(G-uv) + b_{i-2}(G-u-v)+2b_{i-g}(G-C_g)\ \mbox{if}\
g\not\equiv 0  \ (\textnormal{mod} \ 4).
\end{array}
\]

 (2) If $uv$ is a cut edge of $G$, then
 $$b_i(G)=b_i(G-uv) +b_{i-2}(G-u-v).$$
  In particular, if $uv$ is a pendent edge with pendent vertex $u$, then
  $$b_i(G)=b_i(G-u) + b_{i-2}(G-u-v).$$
\end{lem}

\begin{lem}\textnormal{\cite{G1}}\label{lem tree maximal-minimal}
Let $T$ be a tree on $n$ vertices. If $T$ is  different from the
path $P_n$ and the star $S_n$, then $P_n\succ T\succ S_n$.
\end{lem}

 Let $T$ be a tree of order $n\geq 3$, $e=uv$ be a non-pendent edge of $T$. Denote by $T_1$ and $T_2$
 the two components of $T-e$ with $u\in T_1$ and $v\in T_2$. If $T'$
 is the tree obtained from $T$ by contracting the edge $e=uv$ and
 attaching a pendent vertex to the vertex $u$ ($=v$), we say that
 $T'$ is obtained from $T$ by edge-growing transformation (on edge $e=uv$), or e.g.t (on edge $e=uv$) for short \cite{Xu,LGL}.

\begin{lem}\textnormal{\cite{LGL}}\label{lem edge-growing transformation}
If $T'$ is obtained from $T$ by one step of e.g.t, then $T\succ T'$.
\end{lem}

\begin{lem} \textnormal{\cite{SSGL}}\label{lem delete-edge-2}
 Let $G$ be a unicyclic  or bipartite graph, $uv$ be a cut edge of $G$. Then $G\succ G-uv$.
\end{lem}

\begin{lem} \textnormal{\cite{SSGL}}\label{lem delete-vertex}
 Let $G$ be a unicyclic  or bipartite graph, $v$ be a non-isolated vertex in $G$,  and $K_1$ be the trivial graph of order $1$. Then $G\succ (G-v)\cup K_1$.
\end{lem}

\begin{figure}[ht]
\centering
   \setlength{\unitlength}{0.05 mm}%
  \begin{picture}(1197.3, 437.5)(0,0)
  \put(0,0){\includegraphics{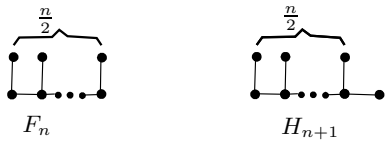}}
  \put(78.18,43.07){\fontsize{8.53}{10.24}\selectfont \makebox(150.0, 60.0)[l]{$F_n$\strut}}
  \put(114.75,329.90){\fontsize{8.53}{10.24}\selectfont \makebox(390.0, 60.0)[l]{$\frac{n}{2}$\strut}}
  \put(765.29,38.20){\fontsize{8.53}{10.24}\selectfont \makebox(270.0, 60.0)[l]{$H_{n+1}$\strut}}
  \put(767.34,335.67){\fontsize{8.53}{10.24}\selectfont \makebox(390.0, 60.0)[l]{$\frac{n}{2}$\strut}}
  \end{picture}%
  \caption{The graphs $F_n$ and $H_{n+1}$.}\label{fig-minimal-tree-degree-3}
\end{figure}

 Let $F_n$ and $H_{n+1}$ be the graphs shown in Figure \ref{fig-minimal-tree-degree-3}.

\begin{lem} \textnormal{\cite{ZL}}\label{lem minimal-tree-degree-3}
 Among conjugated trees of order $n$ with $\Delta\leq 3$, $F_n$ has minimal energy.
\end{lem}


\begin{lem}\label{lem formular for Bn and An}
For $A_n$, $B_n$, $D_n$ and $E_n$, we have that

 (1) $
 b_{2i}(A_n) = b_{2i}(F_n) + b_{2i-2}(F_{n-2})-2 b_{2i-4}(F_{n-4})
            =  b_{2i}(H_{n-1})+2b_{2i-2}(H_{n-3}).
$

(2) $b_{2i}(B_n) = b_{2i}(A_n) + 2b_{2i-6}(F_{n-8}).$

 (3)
\vskip -1cm \hskip -2cm
 \[
\begin{array}{lll}
 b_{2i}(A_n)  =  b_{2i}(A_{n-2})+b_{2i-2}(A_{n-2})+b_{2i-2}(A_{n-4}) \ (n\geq
 8),\\
 b_{2i}(B_n)  =  b_{2i}(B_{n-2}) + b_{2i-2}(B_{n-2})+ b_{2i-2}(B_{n-4}) \ (n\geq
 10),\\
 b_{2i}(D_n)  =  b_{2i}(D_{n-2}) + b_{2i-2}(D_{n-2})+ b_{2i-2}(D_{n-4}) \ (n\geq
 10),\\
 b_{2i}(E_n)  =  b_{2i}(E_{n-2}) + b_{2i-2}(E_{n-2})+ b_{2i-2}(E_{n-4}) \ (n\geq
 12).
\end{array}
\]
\end{lem}

 \pf (1) By Lemma \ref{lem delete-edge} (1) and (2), we can obtain
 that
\[
\begin{array}{lll}
 b_{2i}(A_n) & = & b_{2i}(F_n) + b_{2i-2}(F_{n-2})-2 b_{2i-4}(F_{n-4})\\
             & = & b_{2i}(H_{n-1})+2b_{2i-2}(H_{n-3}).
\end{array}
\]

(2) By Lemma \ref{lem delete-edge} (2), we have
 \[
\begin{array}{lll}
 b_{2i}(B_n) & = & b_{2i}(K_2 \cup A_{n-2}) + b_{2i-2}(H_{n-3})\\
             & = & b_{2i}(K_2 \cup A_{n-2}) + b_{2i-2}(F_{n-4})+
             b_{2i-4}(F_{n-6}),\\
 b_{2i}(A_n) & = & b_{2i}(K_2 \cup A_{n-2}) + b_{2i-2}(A_{n-4})\\
             & = & b_{2i}(K_2 \cup A_{n-2}) + b_{2i-2}(F_{n-4})+
             b_{2i-4}(F_{n-6})-2 b_{2i-6}(F_{n-8}).
\end{array}
\]
Hence
 $ b_{2i}(B_n)  =  b_{2i}(A_n)+2 b_{2i-6}(F_{n-8})$.

 (3) The results directly follow from Lemma \ref{lem delete-edge} (2). \qed

\begin{thm}\label{lem Bn and An}
 $B_n\succ A_n$ for $n\geq 8$.
\end{thm}

 \pf It is easy to obtain that
 \[
 \begin{array}{lll}
 \phi(B_8,x)  =  x^8-8x^6+16x^4-8x^2,\\
 \phi(A_8,x)  =  x^8-8x^6+16x^4-6x^2.
\end{array}
\]
 So $B_8\succ A_8$. Suppose $n>
 8$. By Lemma \ref{lem formular for Bn and An},
$b_{2i}(B_n)-b_{2i}(A_n)=2 b_{2i-6}(F_{n-8})\geq 0$ and $b_{8}(B_n)>
b_{8}(A_n)$. So $B_n\succ A_n$ for
 $n>8$. The proof is thus complete.   \qed

\begin{thm}\label{lem En and Dn}
 $D_8\succ E_8$ and $E_n\succ D_n$ for $n\geq 10$.
\end{thm}

 \pf It is easy to obtain that
\[
\begin{array}{lll}
\phi(E_{8},x) =   x^8-8x^6+14x^4-8x^2+1,\\
\phi(D_{8},x) =   x^8-8x^6+15x^4-8x^2+1,\\
\phi(E_{10},x) =   x^{10}-10x^8+29x^6-31x^4+12x^2-1,\\
 \phi(D_{10},x)  =   x^{10}-10x^8+29x^6-28x^4+10x^2-1,\\
 \phi(E_{12},x)  =   x^{12}-12x^{10}+47x^8-74x^6+51x^4-14x^2+1,\\
 \phi(D_{12},x)  =   x^{12}-12x^{10}+47x^8-72x^6+46x^4-12x^2+1.
\end{array}
\]
 So $D_{8}\succ E_{8}$, $E_{10}\succ D_{10}$ and $E_{12}\succ D_{12}$. Suppose $n > 12$. By Lemma \ref{lem delete-edge},
 \[
\begin{array}{lll}
 b_{2i}(E_n) =  b_{2i}(E_{12}\cup F_{n-12}) + b_{2i-2}(E_{10}\cup F_{n-14}),\\
 b_{2i}(D_n)  =  b_{2i}(D_{12}\cup F_{n-12}) + b_{2i-2}(D_{10}\cup E_{F-14}).
\end{array}
\]
Hence we have $E_n\succ D_n$ for $n>12$. The proof is thus complete.
\qed

\begin{thm}\label{lem Sn and Bn}
 $E(S_8^{4})>E(B_8)$ and $S_n^{\frac{n}{2}}\succ B_n$ for $n\geq 10$.
\end{thm}

 \pf It is easy to obtain that $9.65685\doteq E(S_8^{4})>E(B_8)\doteq
 9.15298$.
 Suppose $n\geq 10$. Since $b_{2i+1}(S_n^{\frac{n}{2}})\geq 0=b_{2i+1}(B_n)$, we
 only need to consider $b_{2i}(S_n^{\frac{n}{2}})$ and $b_{2i}(B_n)$.
 By Lemmas \ref{lem delete-edge} and \ref{lem formular for Bn and An},
\[
 b_{2i}(S_n^{\frac{n}{2}}) =  \left\{\begin{array}{ll} b_{2i}(F_n) + b_{2i-2}(F_{n-4})-2
 b_{2i-\frac{n}{2}}(\frac{n}{2}K_1), &  \mbox{if}\ g\equiv 0
\ (\textnormal{mod} \ 4) \\
  b_{2i}(F_n) + b_{2i-2}(F_{n-4})+2 b_{2i-\frac{n}{2}}(\frac{n}{2}K_1), &  \mbox{if}\ g\not\equiv 0
 \ (\textnormal{mod} \ 4)
\end{array}\right.,
\]
 \[
\begin{array}{lll}
  b_{2i}(B_n) & = &  b_{2i}(F_n) + b_{2i-2}(F_{n-2})-2b_{2i-4}(F_{n-4})+2b_{2i-6}(F_{n-8})\\
              & = &  b_{2i}(F_n) + b_{2i-2}(F_{n-4})+b_{2i-4}(F_{n-6})-b_{2i-4}(F_{n-4})+2b_{2i-6}(F_{n-8})\\
              & = &  b_{2i}(F_n) + b_{2i-2}(F_{n-4})-b_{2i-6}(F_{n-6})+b_{2i-6}(F_{n-8})\\
              & = &  b_{2i}(F_n) + b_{2i-2}(F_{n-4})-b_{2i-8}(F_{n-8})-b_{2i-8}(F_{n-10}).
\end{array}
\]
Since \[b_{2i-\frac{n}{2}}(\frac{n}{2}K_1) =
\left\{\begin{array}{ll} 1, &  \mbox{if} \  2i=\frac{n}{2}\\
 0, &  \mbox{otherwise}\end{array}\right.,
\]
 we have $b_{2i}(S_n^{\frac{n}{2}})\geq b_{2i}(B_n)$, $b_8(S_{10}^5)>
 b_8(B_{10})$, and $b_{10}(S_n^{\frac{n}{2}})>
 b_{10}(B_n)$ when $n\geq 12$. The proof is thus complete.
\qed

  In the following, we will show that $E(A_n)<E(D_n)$ by using the Coulson integral formula method, which had successfully been applied to compare the energy of two given graphs by Huo et al., see \cite{Huo1}-\cite{Huo4}.
  Before proving it, we prepare some results as follows.

\begin{lem}\textnormal{\cite{CDS}}\label{lem characteristic polynomial}
 Let $uv$ be an edge of $G$. Then
 $$\phi(G,x) = \phi(G-uv,x)-\phi(G-u-v,x)-2\sum_{C\in\mathcal
 {C}(uv)}\phi(G-C,x),$$
 where $\mathcal{C}(uv)$ is the set of cycles containing $uv$. In particular, if $uv$ is
 a pendent edge with pendent vertex $v$, then $\phi(G,x) =
 x\phi(G-v,x)-\phi(G-u-v,x)$.
\end{lem}

\begin{lem}\textnormal{\cite{Z2}}\label{lem log inequality}
For any real number $X > -1$, we have $$\frac{X}{1 + X}\leq \log(1 +
X)\leq X.$$
 In particular, $\log(1 + X) < 0$ if and only if $X < 0$.
\end{lem}

\begin{lem}\textnormal{\cite{G2}}\label{lem energy difference}
If $G_1$ and $G_2$ are two graphs with the same number of vertices,
then $$E(G_1)-E(G_2) =\frac{1}{\pi}
\int_{-\infty}^{+\infty}\log\left|\frac{\phi(G_1,ix)}{\phi(G_2,ix)}\right|dx.$$
\end{lem}

From Lemma \ref{lem characteristic polynomial}, we can easily obtain
the following lemma.

\begin{lem}\label{lem characteristic polynomial-An and Dn}
$ \phi(A_n,x) = (x^2-1)\phi(A_{n-2},x)-x^2\phi(A_{n-4},x)$ for
$n\geq 8$, and
 $\phi(D_n,x) = (x^2-1)\phi(D_{n-2},x)-x^2\phi(D_{n-4},x)$ for $n\geq 10$.

\end{lem}

 By some easy calculations, we have $\phi(A_6,x) = x^6-6x^4+6x^2$, $\phi(A_8,x) = x^8-8x^6+16x^4-6x^2$, $\phi(D_6,x) = x^6-6x^4+5x^2-1$ and $\phi(D_8,x) = x^8-8x^6+15x^4-8x^2+1$.
 Now for convenience, we define some notations as follows:
\begin{eqnarray*}
 Y_1(x)=\frac{x^2-1+\sqrt{x^4-6 x^2+1}}{2}, & & Y_2(x)=\frac{x^2-1-\sqrt{x^4-6 x^2+1}}{2}, \\
 Z_1(x)=\frac{-x^2-1+\sqrt{x^4+6 x^2+1}}{2}, & & Z_2(x)=\frac{-x^2-1-\sqrt{x^4+6 x^2+1}}{2}, \\
 A_1(x)=\frac{\phi(A_8,x)-Y_2(x)\phi(A_6,x)}{(Y_1(x))^4-(Y_1(x))^2x^2}, & & A_2(x)=\frac{\phi(A_8,x)-Y_1(x)\phi(A_6,x)}{(Y_2(x))^4-(Y_2(x))^2x^2},\\
  B_1(x)=\frac{\phi(D_8,x)-Y_2(x)\phi(D_6,x)}{(Y_1(x))^4-(Y_1(x))^2x^2},
& &
B_2(x)=\frac{\phi(D_8,x)-Y_1(x)\phi(D_6,x)}{(Y_2(x))^4-(Y_2(x))^2x^2},
\end{eqnarray*}
\begin{eqnarray*}
   f_6(x)=x^6+6x^4+6x^2, & &  f_8(x)=x^8+8x^6+16x^4+6x^2,\\
     g_6(x)=x^6+6x^4+5x^2+1, & & g_8(x)=x^8+8x^6+15x^4+8x^2+1.
\end{eqnarray*}
 It is easy to check that $Y_1(ix)=Z_1(x)$, $Y_2(ix)=Z_2(x)$,
\begin{eqnarray*}
 A_1(ix)=\frac{f_8(x)+Z_2(x)f_6(x)}{(Z_1(x))^4+(Z_1(x))^2x^2}, & & A_2(ix)=\frac{f_8(x)+Z_1(x)f_6(x)}{(Z_2(x))^4+(Z_2(x))^2x^2},\\
  B_1(ix)=\frac{g_8(x)+Z_2(x)g_6(x)}{(Z_1(x))^4+(Z_1(x))^2x^2}, & & B_2(ix)=\frac{g_8(x)+Z_1(x)g_6(x)}{(Z_2(x))^4+(Z_2(x))^2x^2},
\end{eqnarray*}
 $Z_1(x)+Z_2(x)=-x^2-1$ and $Z_1(x)Z_2(x)=-x^2$.
 In addition, for
 $x>0$, $0<\frac{Z_1(x)}{x}<1$; for
 $x<0$, $-1<\frac{Z_1(x)}{x}<0$.

\begin{lem}\label{lem recursive formula of characteristic polynomial-An and Dn}
For $n\geq 6$ and $x\neq 0$, the characteristic polynomials of $A_n$
and $D_n$ have the following forms:
$$\phi(A_n,x) =
 A_1(x)(Y_1(x))^{\frac{n}{2}}+A_2(x)(Y_2(x))^{\frac{n}{2}}$$
and $$ \phi(D_n,x) =
 B_1(x)(Y_1(x))^{\frac{n}{2}}+B_2(x)(Y_2(x))^{\frac{n}{2}}.
$$
\end{lem}

\pf By Lemma \ref{lem characteristic polynomial-An and Dn}, we have
 that $\phi(A_n,x), \phi(D_n,x)$ satisfy the recursive formula
 $f(n,x)=(x^2-1)f(n-2,x)-x^2f(n-4,x)$.
 Therefore, the form of the general solution of the linear homogeneous recursive relation is
 $f(n,x)=C_1(x)(Y_1(x))^{\frac{n}{2}}+C_2(x)(Y_2(x))^{\frac{n}{2}}$.
 By some simple calculations, together with the initial values
 $\phi(A_6,x)$ and $\phi(A_8,x)$ ($\phi(D_6,x)$ and $\phi(D_8,x)$, respectively), we can get that $C_i(x)=A_i(x)$ ($C_i(x)=B_i(x)$, respectively),
 $i=1, 2$.
\qed

\begin{thm}\label{lem An and Dn}
$E(A_n) < E(D_n)$ for $n\geq 6$.
\end{thm}

\pf By Lemma \ref{lem energy difference}, we have
$$E(A_n)-E(D_n) =\frac{1}{\pi}
\int_{-\infty}^{+\infty}\log\left|\frac{\phi(A_n,ix)}{\phi(D_n,ix)}\right|dx.$$
From Lemma \ref{lem recursive formula of characteristic
polynomial-An and Dn}, we know that both $\phi(A_n,ix)$ and
$\phi(D_n,ix)$ are polynomials of $x$ with all real coefficients.
 For convenience, we abbreviate $A_k(ix)$, $B_k(ix)$ and $Z_k(x)$
 to $A_k$, $B_k$ and $C_k$ for $k=1,2$, and abbreviate $f_k(x)$ and $g_k(x)$
 to $f_k$ and $g_k$ for $k=6,8$, respectively. In the
 following, we assume that $x\neq 0$.
 We distinguish two cases in terms of the parity of $n/2$.

{\bf Case 1.} $n=4k$ $(k\geq 2)$. Notice that $Z_1Z_2=-x^2$. When
 $n \rightarrow \infty$,
\begin{eqnarray*}
\frac{\phi(A_n,ix)}{\phi(D_n,ix)} =
\frac{A_1Z_1^{\frac{n}{2}}+A_2Z_2^{\frac{n}{2}}}{B_1Z_1^{\frac{n}{2}}+B_2Z_2^{\frac{n}{2}}}
                                 =  \frac{A_2+A_1(\frac{Z_1}{x})^n}{B_2+B_1(\frac{Z_1}{x})^n}
                                \rightarrow  \frac{A_2}{B_2}.
\end{eqnarray*}
We will show that 
\begin{eqnarray*}
\log\left|\frac{\phi(A_n,ix)}{\phi(D_n,ix)}\right|<\log\left|\frac{A_2(ix)}{B_2(ix)}\right|=\log\frac{A_2}{B_2}.
\end{eqnarray*}
Assume that
\begin{eqnarray*} 2\log\left|\frac{\phi(A_n,ix)}{\phi(D_n,ix)}\right|-2\log\left|\frac{A_2}{B_2}\right|=\log\left(1+\frac{F_1(n,x)}{G_1(n,x)}\right).
\end{eqnarray*}
Then we get that $G_1(n,x)=(A_2(ix) \phi(D_n,ix))^2>0$ and
\begin{eqnarray*}
F_1(n,x) & = &
  B_2^2(A_1Z_1^{\frac{n}{2}}+A_2Z_2^{\frac{n}{2}})^2-A_2^2(B_1Z_1^{\frac{n}{2}}+B_2Z_2^{\frac{n}{2}})^2\\
   & = &
   (A_1^2B_2^2-A_2^2B_1^2)Z_1^n
    +2A_2B_2(A_1B_2-A_2B_1)(Z_1Z_2)^{\frac{n}{2}}\\
   & < &  0.
\end{eqnarray*}
 Since $Z_1^{n}>0$,
$(Z_1Z_2)^{\frac{n}{2}}=(-x^2)^{\frac{n}{2}}>0$, $A_2B_2>0$, and by
some elementary calculations, we have
\begin{eqnarray*}A_1B_2-A_2B_1=\frac{(-3x^8-15x^6-8x^4)\sqrt{x^4+6x^2+1}}{x^6(x^4+6x^2+1)}<0,\end{eqnarray*} and
\begin{eqnarray*}A_1B_2+A_2B_1=\frac{4x^4+17x^6+20x^8+3x^{10}}{x^6(x^4+6x^2+1)}>0.\end{eqnarray*} Hence by Lemma \ref{lem log inequality},
\begin{eqnarray*}
 \frac{1}{\pi} \int_{-\infty}^{+\infty}\log \frac{A_2(ix)}{B_2(ix)}dx
 \leq \frac{1}{\pi} \int_{-\infty}^{+\infty} \left(\frac{A_2(ix)}{B_2(ix)}-1\right)dx\doteq \frac{1}{\pi}
 (-0.8538292323)<0.
\end{eqnarray*}
Therefore,
$$E(A_n)-E(D_n)\leq \frac{1}{\pi} \int_{-\infty}^{+\infty}\log \frac{A_2(ix)}{B_2(ix)}dx<0.$$

 {\bf Case 2.} $n=4k+2$ $(k\geq 1)$. We will show that
 $\log\left|\frac{\phi(A_n,ix)}{\phi(D_n,ix)}\right|$ is monotonically decreasing in $n$.
 Assume that
\begin{eqnarray*}
 2\log\left|\frac{\phi(A_{n+4},ix)}{\phi(D_{n+4},ix)}\right|-2\log\left|\frac{\phi(A_{n},ix)}{\phi(D_{n},ix)}\right|=\log\left(1+\frac{F_2(n,x)}{G_2(n,x)}\right).
\end{eqnarray*}
 Then we can obtain that
 $G_2(n,x)=(\phi(A_n,ix)\phi(D_{n+4},ix))^2>0$ and
\begin{eqnarray*}
 F_2(n,x) & = & (\phi(A_{n+4},ix)\phi(B_n,ix))^2-(\phi(A_n,ix)\phi(D_{n+4},ix))^2\\
          & = &
          (A_1Z_1^{\frac{n}{2}+2}+A_2Z_2^{\frac{n}{2}+2})^2(B_1Z_1^{\frac{n}{2}}+B_2Z_2^{\frac{n}{2}})^2 \\
          & & -  (B_1Z_1^{\frac{n}{2}+2}+B_2Z_2^{\frac{n}{2}+2})^2(A_1Z_1^{\frac{n}{2}}+A_2Z_2^{\frac{n}{2}})^2\\
          & = & -(A_1B_2-A_2B_1)x^n(Z_1^2-Z_2^2)\cdot F_3(n,x),
\end{eqnarray*}
where $
F_3(n,x)=(A_1B_2+A_2B_1)(Z_1Z_2)^{\frac{n}{2}}(Z_1^2+Z_2^2)+2A_1B_1Z_1^{n+2}+2A_2B_2Z_2^{n+2}$.
Since $A_1B_2-A_2B_1<0$, $Z_1^2-Z_2^2<0$ and $x^n>0$, to prove
$F_2(n,x)<0$, it suffice to show that $F_3(n,x)>0$.
By some elementary calculations, we can get that
\begin{eqnarray*}
F_3(n,x)=\frac{2f_8g_8F_4(n,x)+2f_6g_6F_5(n,x)+(f_8g_6+f_6g_8)F_6(n,x)}{x^4+6x^2+1},
\end{eqnarray*}
where
\begin{eqnarray*}
 F_4(n,x) & = &
 Z_1^{n-4}+Z_2^{n-4}-(Z_1Z_2)^{\frac{n}{2}-3}(Z_1^2+Z_2^2)\\
          & = &
          (Z_1^{\frac{n}{2}-1}-Z_2^{\frac{n}{2}-1})(Z_1^{\frac{n}{2}-3}-Z_2^{\frac{n}{2}-3}),\\
F_5(n,x) & = &
 Z_2^2Z_1^{n-4}+Z_1^2Z_2^{n-4}-(Z_1Z_2)^{\frac{n}{2}-2}(Z_1^2+Z_2^2)\\
          & = &
          (Z_1Z_2)^2(Z_1^{\frac{n}{2}-2}-Z_2^{\frac{n}{2}-2})(Z_1^{\frac{n}{2}-4}-Z_2^{\frac{n}{2}-4}),\\
F_6(n,x) & = &
 2Z_2Z_1^{n-4}+2Z_1Z_2^{n-4}-(Z_1+Z_2)(Z_1Z_2)^{\frac{n}{2}-3}(Z_1^2+Z_2^2)\\
          & = &
          Z_1Z_2[(Z_1^{\frac{n}{2}-1}-Z_2^{\frac{n}{2}-1})(Z_1^{\frac{n}{2}-4}-Z_2^{\frac{n}{2}-4})+(Z_1^{\frac{n}{2}-3}-Z_2^{\frac{n}{2}-3})(Z_1^{\frac{n}{2}-2}-Z_2^{\frac{n}{2}-2})].
\end{eqnarray*}
 Notice that $Z_1>0$, $Z_2<0$, we have $Z_1^{k}-Z_2^{k}>0$ when $k$ is
 odd. On the other hand, for $k\geq 1$, we have
\begin{eqnarray*}
Z_1^{2k}-Z_2^{2k}=Z_1^{2k}\left[1-\left(\frac{x}{Z_1}\right)^{4k}\right]<0.
\end{eqnarray*}
 Therefore, we have $F_4(n,x)\geq 0$, $F_5(n,x)\geq 0$,
 $F_6(n,x)>0$, and so $F_3(n,x)>0$ and $F_2(n,x)<0$.
 Hence,
 $$E(A_n)-E(D_n)\leq E(A_6)-E(D_6)\doteq 6.60272-7.20775<0.$$
The proof is thus complete. \qed

\section{Main results}

 Let $G$ be a graph in $\mathbb{U}_{n}$ with vertex set $V$ and the unique cycle $C_g$. We use $M_G$ to denote one arbitrary selected prefect matching of $G$.
 Let $x, y\in V$. Denote by $d_G(x,y)$ ($d_G(x,C_g)$, respectively) the distance between vertex $x$ and $y$ ($C_g$, respectively).
 Define $d=d(G)=\max_{x\in V(G)}\{d_G(x,C_g)\}$, $V_1=V_1(G)=\{x\in V|d_G(x,C_g)=d(G)\}$, and  $t=t(G)=|V_1(G)|$.
 Clearly, the vertices in $V_1$ are pendent vertices when $d(G)\geq 1$.

\begin{lem}\label{lem girth at least 8-d equals 0}
 Let  $n\equiv  0$ \textnormal{(mod $4$)}, $n\geq 8$. Then $C_n\succ B_n$.
\end{lem}

\pf By Lemma \ref{lem delete-edge},
\[
 b_{2i}(C_n) =  \left\{\begin{array}{ll} b_{2i}(P_n) + b_{2i-2}(P_{n-2}), &  \mbox{if}\ 2i\neq n \\
   b_{2i}(P_n) + b_{2i-2}(P_{n-2})-2, &  \mbox{if}\  2i=n
\end{array}\right.,
\]
and
 $$b_{2i}(B_n)=b_{2i}(G_1) + b_{2i-2}(K_2\cup F_{n-4})-2b_{2i-4}(K_2\cup F_{n-6}),$$
 where $G_1$ is the tree of order $n$ obtained by attaching a path with $4$ edges
 to one of the two vertices of degree $2$ of $F_{n-4}$.
 By Lemmas \ref{lem tree maximal-minimal} and \ref{lem
 delete-edge-2}, $P_n\succ G_1$, $P_{n-2}\succ F_{n-2}\succ K_2\cup
 F_{n-4}$. On the other hand, we have $b_{n}(C_n)=b_{n}(B_n)$ and
 $b_{4}(C_n)>b_{4}(B_n)$. Hence $C_n\succ B_n$. \qed

\begin{lem}\label{lem girth at least 8-d equals 1}
 Let $G\in \mathbb{U}_{n}$, $g(G)\equiv  0$ \textnormal{(mod $4$)}, $g(G)\geq 8$, and $d(G)=1$. Then $G\succ B_n$.
\end{lem}

 \pf It is easy to see that $n\geq 10$ and $t(G)$ is even. If $t=g$, then $G=S_n^{\frac{n}{2}}$, and so $G\succ
 B_n$ by Lemma \ref{lem Sn and Bn}. So in the following, we suppose $2\leq t\leq g-2$.
 Let $C_g=y_1y_2\ldots y_g$, and $x_{k_1}y_{k_1}$, $x_{k_2}y_{k_2}$,
 \ldots, $x_{k_t}y_{k_t}$ ($1=k_1< k_2<\cdots < k_t$) be all the edges outside $C_G$. Then there
 must exist an index $k_i$ such that $k_{i}+1, k_{i}+2$ (mod $g$)
  are not in the set $\{k_1, k_2, \ldots, k_t\}$, that is $y_{k_{i}+1}y_{k_{i}+2}\in M_G$.
 Without loss of generality, we assume that $k_i=1$. Since $G$ has a perfect matching, we have that $4\leq k_2 < k_3 < \cdots < k_{t}\leq g$, and $k_2-3$ is odd, $k_{i+1}-k_i$ is odd for $2\leq i\leq t-1$, and $g-k_{t}$ is even.

 By Lemma \ref{lem delete-edge}, we have $ b_{2i}(G) = b_{2i}(G-x_1) + b_{2i-2}(G-x_1-y_1)$, and
\[
b_{2i}(G-x_1)  =  \left\{\begin{array}{ll}b_{2i}(G-x_1-y_1y_2)+ b_{2i-2}(G-x_1-y_1-y_2), &  \mbox{if}\ 2i\neq g \\
b_{2i}(G-x_1-y_1y_2)+ b_{2i-2}(G-x_1-y_1-y_2)-2, &  \mbox{if}\ 2i=g
\end{array}\right..
\]
  Note that $G-x_1-y_1$ is a conjugated tree of order $n-2$ with $\Delta\leq 3$, by Lemma \ref{lem minimal-tree-degree-3}, $G-x_1-y_1\succeq F_{n-2}$.

 Denote $T_1=G-x_1-y_1y_2$ and $T_2=G-x_1-y_1-y_2$.  Notice that if $k_{i+1}-k_i=2k+1>1$,
 then from $T_1$, we can obtain a different tree $T_2$ with $\Delta \leq 3$ by carrying out $k$ steps of
 e.g.t. Therefore we can finally get $H_{n-1}$ from $T_1$ by carrying out
 e.g.t repeatedly, if necessary. By Lemma \ref{lem edge-growing
 transformation}, we have $T_1\succeq H_{n-1}$. Similarly, we can
 obtain that $T_2\succeq H_{n-3}$.
 By Lemma \ref{lem formular for Bn and An}, we have $$b_{2i}(B_n)= b_{2i}(H_{n-1}) +
 2b_{2i-2}(H_{n-3})+2 b_{2i-6}(F_{n-8}).$$
 Hence for $2i\neq g$, we have
 \[
\begin{array}{lll}
  b_{2i}(G) & \geq &  b_{2i}(H_{n-1}) + b_{2i-2}(F_{n-2})+b_{2i-2}(H_{n-3})\\
              & = & b_{2i}(H_{n-1}) + 2b_{2i-2}(H_{n-3})+b_{2i-4}(F_{n-4})\\
              & = & b_{2i}(B_n)+b_{2i-4}(F_{n-4})-2b_{2i-6}(F_{n-8})\\
              & = &
              b_{2i}(B_n)+b_{2i-4}(F_{n-6})+b_{2i-8}(F_{n-8})+b_{2i-8}(F_{n-10})\\
              & \geq & b_{2i}(B_n),
\end{array}
\]
and
$$
  b_{g}(G)  \geq  b_{g}(B_n)+b_{g-4}(F_{n-6})+b_{g-8}(F_{n-8})+b_{g-8}(F_{n-10})-2\geq  b_{g}(B_n),$$
 since $g\leq n-2$, $b_{g-4}(F_{n-6})\geq 1$ and $b_{g-8}(F_{n-8})\geq 1$.
On the other hand, it is obvious that $b_{4}(G) >  b_{4}(B_n)$. Thus
$G\succ B_n$.   \qed

\begin{lem}\label{lem girth at least 8-d equals 2}
 Let $G\in \mathbb{U}_{n}$, $g(G)\equiv  0$ \textnormal{(mod $4$)}, $g(G)\geq 8$, and $d(G)=2$. Then $G\succ B_n$.
\end{lem}

 \pf  Let $C_g=z_1z_2\ldots z_g$. We apply induction on $t$.
 Suppose $t=1$. Assume that $d_G(x_1,z_1)=2$ and $x_1y_1z_1$ be a path with length $2$ in $G$. Then $d_G(y_1)=2$ and $x_1y_1\in M_G$.
 Since $G$ has a perfect matching, either $z_1z_2\in M_G$ or $z_1z_g\in M_G$. We may assume that $z_1z_2\in M_G$.  Let $E_1=\{z_{i_1}y_{i_1}, z_{i_2}y_{i_2}, \ldots, z_{i_k}y_{i_k}\}$ ($i_1<i_2<\cdots <i_k$) be the
 set of all edges in $M_G\setminus x_1y_1$ outside $C_G$. 
 If $E_1$ is not empty, then we have that $3\leq i_1 < i_2 < \cdots < i_{k}\leq g$, and $i_1-3$ is even, $i_{j+1}-i_j$ is odd for $1\leq j\leq k-1$, and $g-i_{k}$ is even.
By Lemma \ref{lem delete-edge}, we have
\[
\begin{array}{lll}
 b_{2i}(G) & = &  b_{2i}(G-x_1) + b_{2i-2}(G-x_1-y_1)\\
           & = &  b_{2i}(G-x_1-y_1) + b_{2i-2}(G-x_1-y_1-z_1) +b_{2i-2}(G-x_1-y_1).
\end{array}
\]
 Denote $G_1=G-x_1-y_1$ and $G_2=G-x_1-y_1-z_1$. Then $G_1\in\mathbb{U}_{n-2}$ with $d(G_1)\leq
 1$. By Lemmas \ref{lem girth at least 8-d equals 0} and \ref{lem girth at least 8-d equals 1}, we have $G_1\succ B_{n-2}$.
 By an argument similar to the proof in Lemma \ref{lem girth at least 8-d equals
 1}, we can obtain $G_2\succeq H_{n-3}$. Therefore by Lemmas \ref{lem
 delete-edge} and \ref{lem formular for Bn and An},
\[
\begin{array}{lll}
 b_{2i}(G) & \geq &  b_{2i}(B_{n-2}) + b_{2i-2}(B_{n-2})+ b_{2i-2}(H_{n-3}) \\
           & = & b_{2i}(B_{n-2}) + b_{2i-2}(B_{n-2})+ b_{2i-2}(H_{n-5})+ 2b_{2i-4}(H_{n-7}) + 2b_{2i-6}(F_{n-8})\\
           & = & b_{2i}(B_{n-2}) + b_{2i-2}(B_{n-2})+ b_{2i-2}(A_{n-4}) + 2b_{2i-6}(F_{n-8}) \\
           & = & b_{2i}(B_{n-2}) + b_{2i-2}(B_{n-2})+ b_{2i-2}(B_{n-4}) + 2b_{2i-6}(F_{n-8})-2b_{2i-8}(F_{n-12}) \\
           & = & b_{2i}(B_{n-2}) + b_{2i-2}(B_{n-2})+ b_{2i-2}(B_{n-4}) + 2b_{2i-6}(F_{n-10})+2b_{2i-8}(F_{n-10}) \\
           & = & b_{2i}(B_{n}) + 2b_{2i-6}(F_{n-10})+2b_{2i-8}(F_{n-10}) \\
           & \geq & b_{2i}(B_{n}).
\end{array}
\]
 Since $G_1\succ B_{n-2}$, there exist $i$ such that
 $b_{2i}(G)>b_{2i}(B_{n})$. Hence $G\succ B_n$ for $t=1$.

 Assume now that $t\geq 2$ and the assertion holds for smaller values of $t$.
 Let $V_1(G)=\{x_1, x_2, \ldots, x_t\}$, and for each $i\in \{1, 2, \ldots,t \}$, $x_iy_iz_{k_i}$ ($1=k_1<k_2<\cdots<k_t\leq g$) be a path with length $2$ in $G$. Then $d_G(y_i)=2$ and $x_iy_i\in M_G$.
 For convenience, let $k_{t+1}=k_1$. We consider the following two cases:

 {\bf Case 1.} There exist two indices $k_i$ and $k_{i+1}$ such that $z_{k_i}z_{k_{i+1}}$ is an edge on $C_g$.
 Without loss of generality, we assume that $k_i=1, k_{i+1}=2$.
 Now let $G'$ be the graph obtained from $G$ by deleting the edge
 $z_2y_2$ and adding one new edge $y_2y_1$. Then $G'\in \mathbb{U}_{n}$,
 and $d(G')=3$.

 {\bf Claim 1.} $G\succ G'$.

 \pf By Lemma \ref{lem delete-edge}, we have
\[
\begin{array}{lll}
  b_{2i}(G) & = &  b_{2i}(G-x_2) + b_{2i-2}(G-\{x_2,y_2\}),\\
  b_{2i}(G') & = &  b_{2i}(G'-x_2) + b_{2i-2}(G'-\{x_2,y_2\}),
\end{array}
\]
and
\[
\begin{array}{lll}
b_{2i}(G-x_2)  & = & b_{2i}(G-\{x_2,x_1\})+ b_{2i-2}(G-\{x_2,x_1,y_1\})\\
               & = &
                b_{2i}(G-\{x_2,x_1,y_2\})+b_{2i-2}(G-\{x_2,x_1,y_2,z_2\})\\
               &  &
               +b_{2i-2}(G-\{x_2,x_1,y_1,y_2\})+b_{2i-4}(G-\{x_2,x_1,y_1,y_2,z_2\}),\\
b_{2i}(G'-x_2)  & = & b_{2i}(G'-\{x_2,x_1\})+ b_{2i-2}(G'-\{x_2,x_1,y_1\})\\
               & = &
                b_{2i}(G'-\{x_2,x_1,y_2\})+
                2b_{2i-2}(G'-\{x_2,x_1,y_2,y_1\}).
\end{array}
\]
  Denote $G_3=G-\{x_2,x_1,y_2,z_2\}$, $G_4=G-\{x_2,x_1,y_1,y_2,z_2\}$
 and $G_5=G-\{x_2,x_1,y_2,y_1\}$. Since $G-\{x_2,y_2\}=G'-\{x_2,y_2\}$,
 $G-\{x_2,x_1,y_2\}=G'-\{x_2,x_1,y_2\}$, and $G-\{x_2,x_1,y_1,y_2\}=G'-\{x_2,x_1,y_2,y_1\}$, we have
$$
b_{2i}(G)- b_{2i}(G') =
        b_{2i-2}(G_3)+b_{2i-4}(G_4) -b_{2i-2}(G_5).
$$
Furthermore,
\[
\begin{array}{lll}
b_{2i-2}(G_5) & = & b_{2i-2}(G_5-\{z_1z_2\})+b_{2i-4}(G_5-\{z_1,z_2\})-2b_{2i-g}(G_5-C_g)\\
              & = &
              b_{2i-2}(G_5-\{z_1z_2,z_2z_3\})+b_{2i-4}(G_5-\{z_1z_2\}-\{z_2,z_3\})\\
            &  &     +b_{2i-4}(G_5-\{z_1,z_2\})-2b_{2i-g}(G_5-C_g)\\
           & = &  b_{2i-2}(G_5-\{z_2\})+b_{2i-4}(G_5-\{z_2,z_3\})\\
             &  &
             +b_{2i-4}(G_5-\{z_1,z_2\})-2b_{2i-g}(G_5-C_g),\\
b_{2i-2}(G_3) & = & b_{2i-2}(G_3-\{y_1\})+b_{2i-4}(G_3-\{y_1,z_1\})\\
              & = & b_{2i-2}(G_5-\{z_2\})+b_{2i-4}(G_5-\{z_1,z_2\}).
\end{array}
\]
 Note that $G_5-\{z_2,z_3\}=G_4-z_3$. Hence $G_4\succ [G_5-\{z_2,z_3\}]\cup
 K_1$ by Lemma \ref{lem delete-vertex}.
 Therefore $$b_{2i}(G)- b_{2i}(G') = b_{2i-4}(G_4) - b_{2i-4}(G_5-\{z_2,z_3\})+2b_{2i-g}(G_5-C_g)\geq 0,$$
 and $$b_{6}(G)- b_{6}(G') = b_{2}(G_4)-b_{2}(G_5-\{z_2,z_3\})>0.$$
Thus the result $G\succ G'$ holds. \qed

 Now it suffices to prove that $G'\succeq B_n$. By Lemma \ref{lem delete-edge}, we have
\[
\begin{array}{lll}
  b_{2i}(G') & = &  b_{2i}(G'-x_2) + b_{2i-2}(G'-\{x_2,y_2\})\\
             & = &  b_{2i}(G'-\{x_2,y_2\}) + b_{2i-2}(G'-\{x_2,y_2,y_1\})+ b_{2i-2}(G'-\{x_2,y_2\})\\
             & = &  b_{2i}(G'-\{x_2,y_2\}) + b_{2i-2}(G'-\{x_2,y_2,y_1,x_1\})+
             b_{2i-2}(G'-\{x_2,y_2\}).
\end{array}
\]
 Denote $G_6=G'-\{x_2,y_2\}$ and $G_7=G'-\{x_2,y_2,y_1,x_1\}$. Then it is easy to see that $G_6\in\mathbb{U}_{n-2}$, $G_7\in\mathbb{U}_{n-4}$,
 and for $i=6,7$, we have $d(G_i)\leq 1$ or $d(G_i)=2$ and $t(G_i)<t(G)$. Therefore by Lemmas \ref{lem girth at least 8-d equals 0}, \ref{lem girth at least 8-d equals 1} and the induction
 hypothesis, $G_6\succ B_{n-2}$, $G_7\succ B_{n-4}$ and so we have $G'\succ B_n$.

 {\bf Case 2.} For the indices $1=k_1<k_2<\ldots<k_t\leq g$, we have $k_{i+1}-k_i\geq 2$ for $1\leq i\leq t-1$ and $k_t\leq g-2$.
  Since $G$ has a perfect matching, we may assume that $z_1z_2\in M_G$, and $z_{k_i}w_{k_i}\in M_G$ for $2\leq i\leq t$, where $w_{k_i}\in\{z_{k_i-1},z_{k_i+1}\}$.
  Similarly, we have
  $$b_{2i}(G) =  b_{2i}(G-x_1-y_1) + b_{2i-2}(G-x_1-y_1-z_1) + b_{2i-2}(G-x_1-y_1),$$
  and
  $$G-x_1-y_1\succ B_{n-2}.$$
 Denote $T=G-x_1-y_1-z_1$.
 Let $T_1$ be the tree obtained from $T$ by
 deleting the edge $x_{i}y_{i}$ and adding one new edge $x_iw_{{k_i}}$ ($2\leq i\leq t$), we say that $T_1$ is obtained from $T$ by Operation I.

 {\bf Claim 2.} $T\succ T_1$.

 \pf By Lemma \ref{lem delete-edge}, we have
\[
\begin{array}{lll}
 b_{2i}(T)  & = & b_{2i}(T-x_i)+ b_{2i-2}(T-\{x_i,y_i\}), \\
 b_{2i}(T_1)  & = & b_{2i}(T_1-x_i)+b_{2i-2}(T_1-\{x_i,w_{{k_i}}\}).
\end{array}
\]
 Note that $T-x_i=T_1-x_i$, and $T_1-\{x_i,w_{{k_i}}\}$ is isomorphic to a proper subgraph of
 $T-\{x_i,y_i\}$, then $T-\{x_i,y_i\}\succ T_1-\{x_i,w_{{k_i}}\}$ by Lemma \ref{lem
 delete-edge-2}.  Therefore $T\succ T_1$. \qed

 Let $T'$ be the tree obtained from $T$ by deleting $t-1$ edges $x_{2}y_{2}, x_{3}y_{3},\ldots, x_ty_t$ and adding $t$ new edges $x_{2}w_{{k_2}}, x_{3}w_{{k_3}},\ldots, x_tw_{{k_t}}$.
 Then from $T$ we can obtain $T'$ by applying Operation I $t-1$ times.
 By Claim 2, we have $T\succ T'$.
 Clearly, $T'$ is a tree with $\Delta\leq 3$. Now we can assume that $z_{j_1}, z_{j_2},\ldots, z_{j_{l}}$ ($2<j_1 < j_2 < \cdots < j_{l}\leq g$) are all vertices with degree $3$ in $T'$. Then we have
 $j_1-2$ is odd, $j_{i+1}-j_i$ is odd for $1\leq i\leq l-1$, and $g-j_{l}$ is even.
 By an argument similar to the proof in Lemma \ref{lem girth at least 8-d equals 1}, we can obtain $T'\succeq H_{n-3}$.
 Therefore we have $T\succ H_{n-3}$ and so similar to the above case $t=1$, we
 can finally obtain $G\succ B_n$.

 The proof is thus complete.\qed

 Let $G$ be a graph in $\mathbb{U}_{n}$ with $g(G)\equiv  0$ \textnormal{(mod $4$)}, $d(G)\geq 3$.
 Suppose $C_g=z_1z_2\ldots z_g$, $x_1,y_1\in V_1(G)$, and $x_1x_2 x_3\ldots x_dz_1$, $y_1y_2x_3 \ldots x_dz_1$ be two paths with length $d$ in
 $G$. For convenience, denote $x_{d+1}=z_1$.
 If $G'$ is the graph obtained from $G$ by deleting two edges $x_3y_2$, $y_2y_1$ and adding two new edges $y_1x_1$ and
 $y_2x_2$, then we say that $G'$ is obtained from $G$ by Operation
 II.  Clearly, $G'\in\mathbb{U}_{n}$.

\begin{lem}\label{lem Operation 2}
 Let $G$ be defined as above. If $G'$ is obtained from $G$ by Operation II, then $G\succ G'$.
\end{lem}

\pf By Lemma \ref{lem delete-edge}, we have
\[
\begin{array}{lll}
  b_{2i}(G) & = &  b_{2i}(G-y_1) + b_{2i-2}(G-\{y_1,y_2\})\\
            & = & b_{2i}(G-\{y_1,y_2\})+ b_{2i-2}(G-\{y_1,y_2,x_3\}) +
            b_{2i-2}(G-\{y_1,y_2\}),\\
  b_{2i}(G') & = &  b_{2i}(G'-y_1) + b_{2i-2}(G'-\{y_1,x_1\})\\
             & = &   b_{2i}(G'-\{y_1,y_2\})+ b_{2i-2}(G'-\{y_1,y_2,x_2) +
             b_{2i-2}(G'-\{y_1,x_1\}).
\end{array}
\]
 It is easy to see that $G'-\{y_1,y_2\}=G-\{y_1,y_2\}\cong
 G'-\{y_1,x_1\}$. Let $G_1=G-\{x_1,x_2,y_1,y_2\}$. Then
\[
\begin{array}{lll}
 b_{2i-2}(G-\{y_1,y_2,x_3\})  & = &  b_{2i-2}((G_1-x_3)\cup K_2)\\
                             & = &  b_{2i-2}(G_1-x_3)+b_{2i-4}(G_1-x_3),\\
 b_{2i-2}(G'-\{y_1,y_2,x_2\})  & = & b_{2i-2}(G_1)\\
                             & = &  b_{2i-2}(G_1-x_3)+b_{2i-4}(G_1-\{x_3,x_4\}).\\
\end{array}
\]
 Since $G_1-x_3\succ (G_1-\{x_3,x_4\})\cup K_1$ by Lemma \ref{lem
 delete-vertex}, we have $b_{2i}(G)\geq b_{2i}(G')$ and $b_{6}(G)>
 b_{6}(G')$. Thus the result $G\succ G'$ holds. \qed

 Let $G$ be a graph in $\mathbb{U}_{n}$ with $g(G)\equiv  0$ \textnormal{(mod $4$)}, $d(G)=3$.
 Suppose $C_g=z_1z_2\ldots z_g$, $x_1,y_1\in V_1(G)$, and $x_1x_2 x_3z_1$, $y_1y_2y_3z_i$ ($2\leq i\leq g$) be two paths with length $3$ in
 $G$, where $d_G(x_3)=d_G(y_3)=2$.
 If $G'$ is the graph obtained from $G$ by deleting two edges $y_1y_2$, $y_2y_3$ and adding two new edges $y_1x_1$ and
 $y_2x_2$, then we say that $G'$ is obtained from $G$ by Operation
 III.  Clearly, $G'\in\mathbb{U}_{n}$.

\begin{lem}\label{lem Operation 3}
 Let $G$ be defined as above. If $G'$ is obtained from $G$ by Operation III, then $G\succ G'$.
\end{lem}

\pf The proof is similar to that of Lemma \ref{lem Operation 2}.
\qed

\begin{thm}\label{thm girth at least 8}
 Let $G\in \mathbb{U}_{n}$, $g(G)\equiv  0$ \textnormal{(mod $4$)}, $g(G)\geq 8$. Then $G\succ B_n$.
\end{thm}

 \pf We apply induction on $d$. As the case $d\leq 2$ was proved by Lemmas \ref{lem girth at least 8-d equals 0}, \ref{lem girth at least 8-d equals 1}, and \ref{lem girth at least 8-d equals 2},
 we now suppose that $d\geq 3$ and the assertion holds for smaller values of $d$.

  Let $C_g=z_1z_2\ldots z_g$. Assume that $d_G(x_1,z_1)=d$ and $x_1x_2 \ldots x_dz_1$ be a path with length $d$ in $G$. Then $d_G(x_2)=2$ and $x_1x_2\in M_G$.
  For convenience, denote $x_{d+1}=z_1$.
  Let $G_1=G-\{x_1,x_2\}$ and $G_2=G-\{x_1,x_2,x_3\}$. By Lemma \ref{lem delete-edge}, we have
  \begin{eqnarray}b_{2i}(G) =  b_{2i}(G_1) + b_{2i-2}(G_1)+ b_{2i-2}(G_2).\label{coefficient formular}
  \end{eqnarray}
 Note that
\begin{eqnarray} b_{2i}(B_n)  = b_{2i}(B_{n-2})+b_{2i-2}(B_{n-2})+b_{2i-2}(B_{n-4}) \ (n\geq 10).\label{coefficient formular-Bn} \end{eqnarray}

  Now we prove the result for the given $d$ by induction on $t$.

  {\bf Case 1.} Suppose $t=1$.

   {\bf Subcase 1.1.} $d_G(x_3)=d_G(x_4)=2$.

   Note that $x_3x_4\in M_G$. It is easy to see that $G_1\in \mathbb{U}_{n-2}$,
  $G_2-x_4\in \mathbb{U}_{n-4}$, $d(G_1)<d$ and $d(G_2-x_4)<d$. So by the induction hypothesis,
  $G_1\succeq
  B_{n-2}$, $G_2-x_4\succeq B_{n-4}$. By Lemma \ref{lem delete-vertex}, we have $G_2\succ (G_2-x_4)\cup K_1$, and so $G_2\succ B_{n-4}\cup K_1$.
  It follows from Eqs. \eqref{coefficient formular} and \eqref{coefficient
  formular-Bn} that $G\succ B_n$.

   {\bf Subcase 1.2.} $d_G(x_3)=2, d_G(x_4)=3$ and $d>3$.

  Note that $x_3x_4\in M_G$. Suppose $y_2\not\in\{x_3,x_5\}$ is a neighbor of $x_4$, and $y_1y_2\in M_G$. Since $t=1$, we
  have $d_G(y_1)=1$ and $d_G(y_2)=2$. Then we have $G_2-y_1\in \mathbb{U}_{n-4}$ since $(M_G\setminus(\{x_1x_2,x_3x_4,y_1y_2\})\cup\{y_2x_4\}$ is a perfect matching of $G_2-y_{1}$.
   Therefore similar to Subcase 1.1, we have $G_1\succeq B_{n-2}$ and $G_2\succ B_{n-4}\cup
   K_1$, and so $G\succ B_n$.

  {\bf Subcase 1.3.} $d_G(x_3)=2, d_G(x_4)=3$ and $d=3$, i.e., $x_{4}=z_1$.

  Note that $x_3x_4\in M_G$. Since $G$ has a perfect matching, $g(G)\equiv  0$ \textnormal{(mod $4$)}, there exist $k\geq 1$ pendent edges
 $z_{i_1}y_{i_1},\ldots,z_{i_k}y_{i_k}$ ($2\leq i_1<\ldots<i_k\leq g$) such that $k$ is odd and $z_{i_k}y_{i_k}\in M_G$.
 Then we have $G_2-y_{i_1}\in \mathbb{U}_{n-4}$, since
 $$(M_G\setminus\{x_1x_2,x_3x_4,z_{i_1}y_{i_1},z_2z_3,z_4z_5,\ldots,z_{i_1-2}z_{i_1-1}\})\cup\{z_1z_2,z_3z_4,\ldots,z_{i_1-1}z_{i_1}\}$$
 is a perfect matching of $G_2-y_{i_1}$.
 Therefore similarly, we have $G_1\succeq B_{n-2}$ and $G_2\succ B_{n-4}\cup K_1$, and so $G\succ B_n$.

 {\bf Subcase 1.4.} $d_G(x_3)=3$.

 Since $t=1$, $x_3x_4\not\in M_G$, we assume that $x_3y_1\in M_G$.
 Then $d_G(y_1)=1$, $d_{G_2}(y_1)=0$ and $G_2-y_1\in \mathbb{U}_{n-4}$. Therefore similarly, we have $G_1\succ B_{n-2}$ and $G_2\succeq B_{n-4}\cup K_1$, and so $G\succ B_n$.

 {\bf Case 2.} Assume now that $t\geq 2$ and the assertion holds for smaller values of $t$.
 Note that $G_1\in \mathbb{U}_{n-2}$ with $d(G_1)=d$ and $t(G_1)<t$. By
 the induction hypothesis, $G_1\succeq B_{n-2}$.

  {\bf Subcase 2.1.} $d_G(x_3)=d_G(x_4)=2$.

 The proof is similar to that of Subcases 1.1.

  {\bf Subcase 2.2.} $d_G(x_3)=3$ and $x_3x_4\not\in M_G$.

  The proof is similar to that of Subcases 1.4.

  {\bf Subcase 2.3.} $d_G(x_3)=2$, $d_G(x_4)=3$, and $d>3$.

 Suppose $y_2\not\in\{x_3,x_5\}$ is a neighbor of $x_4$, and $y_1y_2\in M_G$. If $d_G(y_2)=3$, let $y_3\not\in\{x_4,y_1\}$ be a neighbor of $y_2$, and $y_4y_3\in M_G$.
 Then we have $d_G(y_1)=d_G(y_4)=1$, $d_G(y_3)=2$, and $G_2-y_1\in \mathbb{U}_{n-4}$, since $(M_G\setminus(\{x_1x_2,x_3x_4,y_1y_2\})\cup\{y_2x_4\}$ is a perfect matching of $G_2-y_{1}$.
 Therefore similarly, we have $G_2\succ B_{n-4}\cup K_1$, and so $G\succ B_n$.

  {\bf Subcase 2.4.} $d_G(x_3)=3$ and $x_3x_4\in M_G$.

 Suppose $y_2\not\in\{x_2,x_4\}$ is a neighbor of $x_3$, and $y_1y_2\in
 M_G$. Then we have $d_G(y_1)=1$ and $d_G(y_2)=2$. Let $G'$ be the graph obtained from $G$ by Operation II. It follows from
 Lemma \ref{lem Operation 2} that $G\succ G'$.
 Similarly, we have  $$b_{2i}(G') =  b_{2i}(G'-\{y_1,x_1\}) + b_{2i-2}(G'-\{y_1,x_1\})+ b_{2i-2}(G'-\{y_1,x_1,x_2,y_2\}),$$
 and $G'-\{y_1,x_1\}\succ B_{n-2}$, and $G'-\{y_1,x_1,x_2,y_2\}\succeq B_{n-4}$. Therefore $G'\succ B_n$.

 {\bf Subcase 2.5.} $d_G(x_3)=2$ and $d=3$.

 Now $x_4=z_1$. Since $t\geq 2$, suppose that $y_1\in V_1(G)$ and $y_1y_2y_3z_i$($i\neq 2$) be a path with length $3$ in
 $G$. Then by the above subcases, we may assume that $d_G(y_3)=2$.
 Let $G'$ be the graph obtained from $G$ by Operation III, then we have $G\succ G'$ by Lemma \ref{lem Operation 3}. And similar to the Subcase 2.4, we
 have $G'\succ B_n$.

 The proof is thus complete. \qed

\begin{figure}[ht]
\centering
     \setlength{\unitlength}{0.05 mm}%
  \begin{picture}(1831.4, 291.0)(0,0)
  \put(0,0){\includegraphics{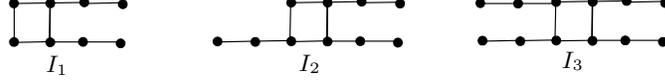}}
  \put(135.76,38.20){\fontsize{8.53}{10.24}\selectfont \makebox(150.0, 60.0)[l]{$I_1$\strut}}
  \put(808.17,41.93){\fontsize{8.53}{10.24}\selectfont \makebox(150.0, 60.0)[l]{$I_2$\strut}}
  \put(1507.92,50.63){\fontsize{8.53}{10.24}\selectfont \makebox(150.0, 60.0)[l]{$I_3$\strut}}
  \end{picture}%
  \caption{The graphs considered in Theorem \ref{thm girth 4-A_n}.}\label{fig-lemma-An}
\end{figure}

\begin{thm}\label{thm girth 4-A_n}
 Let $G\in \mathbb{U}_{n}$, $g(G)=4$, $G\not\cong A_n$. If there are just two edges of $M_G$ in $C_4$, then $G\succ A_n$.
\end{thm}

 \pf We apply induction on $d$. Suppose $G\not\cong A_n$ is a graph in $\mathbb{U}_{n}$ with $g=4$, and there are just two edges of $M_G$ in $C_4$.
 For $d\leq 1$, there is nothing to
 prove. Suppose $d=2$. Then $G$ is isomorphic to one of the following
 graphs $B_8$, $I_1$, $I_2$ and $I_3$, as shown in Figure \ref{fig-lemma-An}. By Lemma \ref{lem Bn and An}, we have $B_8\succ A_8$. It is easy to obtain that
\[
\begin{array}{lll}
\phi(I_{1},x) & = &  x^{8}-8x^6+16x^4-9x^2,\\
 \phi(I_{2},x) & = &  x^{10}-10x^8+30x^6-34x^4+12x^2,\\
 \phi(I_{3},x) & = &  x^{12}-12x^{10}+48x^8-84x^6+64x^4-16x^2,\\
  \phi(A_8,x) & = & x^8-8x^6+16x^4-6x^2,\\
 \phi(A_{10},x) & = &   x^{10}-10x^8+30x^6-28x^4+6x^2,\\
 \phi(A_{12},x) & = &  x^{12}-12x^{10}+48x^8-74x^6+40x^4-6x^2.
 \end{array}
\]
 Hence $I_{1}\succ A_{8}$, $I_{2}\succ A_{10}$, and $I_{3}\succ A_{12}$.
 Now suppose that $d\geq 3$ and the assertion holds for smaller values of
 $d$. By an argument similar to the proof of Theorem \ref{thm girth at least
 8}, we can obtain $G\succ A_n$.  \qed

\begin{figure}[ht]
\centering
    \setlength{\unitlength}{0.05 mm}%
  \begin{picture}(2482.2, 296.0)(0,0)
  \put(0,0){\includegraphics{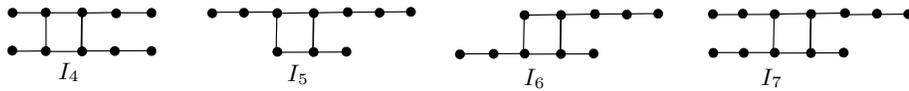}}
  \put(171.40,54.93){\fontsize{8.53}{10.24}\selectfont \makebox(150.0, 60.0)[l]{$I_4$\strut}}
  \put(780.67,48.28){\fontsize{8.53}{10.24}\selectfont \makebox(150.0, 60.0)[l]{$I_5$\strut}}
  \put(1409.27,38.20){\fontsize{8.53}{10.24}\selectfont \makebox(150.0, 60.0)[l]{$I_6$\strut}}
  \put(2040.66,39.44){\fontsize{8.53}{10.24}\selectfont \makebox(150.0, 60.0)[l]{$I_7$\strut}}
  \end{picture}%
  \caption{The graphs considered in Theorem \ref{thm girth 4-D_n}.}\label{fig-lemma-Dn}
\end{figure}

\begin{thm}\label{thm girth 4-D_n}
 Let $G\in \mathbb{U}_{n}$, $g(G)=4$, $G\not\cong D_n$. If there is just one edge of $M_G$ in $C_4$, then $G\succ D_n$.
\end{thm}

 \pf We apply induction on $d$. For $d\leq 1$, there is nothing to prove. Suppose $d=2$. Then $G$ is isomorphic to $I_4$, as shown in Figure \ref{fig-lemma-Dn}.
 It is easy to obtain that
\[
\begin{array}{lll}
\phi(I_{4},x) & = &  x^{10}-10x^8+29x^6-32x^4+12x^2-1,\\
 \phi(D_{10},x) & = &  x^{10}-10x^8+29x^6-28x^4+10x^2-1.
\end{array}
\]
Hence $I_4\succ D_{10}$.

 Now suppose that $d\geq 3$ and the assertion holds for smaller values of
 $d$.
 We use the same notations as in Theorem \ref{thm girth at least 8}. Then the proof is similar to that of Theorem \ref{thm girth at least 8}.
 We can divide two cases $x_dz_1\in M_g$ and $x_dz_1\not\in M_g$ to proceed. The difference is that we need to prove the result $G\succ D_n$ for the case: $d=3$, $t=1$, $d_G(x_3)=2$ and $z_2z_3\in
 M_4$. Since $G\not\cong D_n$, $G$ is isomorphic to one of the following
 graphs $I_5$, $I_6$ and $I_7$, as shown in Figure \ref{fig-lemma-Dn}.
 It is easy to obtain that
\[
\begin{array}{lll}
\phi(I_{5},x) & = &  x^{10}-10x^8+30x^6-33x^4+11x^2-1,\\
\phi(I_{6},x) & = &  x^{10}-10x^8+30x^6-33x^4+12x^2-1,\\
\phi(I_{7},x) & = &  x^{12}-12x^{10}+48x^8-83x^6+62x^4-16x^2+1,\\
 \phi(D_{12},x) & = &  x^{12}-12x^{10}+47x^8-72x^6+46x^4-12x^2+1.
\end{array}
\]
Hence $I_5\succ D_{10}$, $I_6\succ D_{10}$ and $I_7\succ D_{12}$.
\qed

\begin{thm}\label{thm girth 4-E_n}
 Let $G\in \mathbb{U}_{n}$, $g(G)=4$, $G\not\cong E_n$. If there are no edges of $M_G$ in $C_4$, then $G\succ E_n$.
\end{thm}

 \pf  The proof is similar to that of Theorem \ref{thm girth at least 8}.   \qed

\vskip 0.5cm
\noindent {\bf Proof of Theorem \ref{main thm}.}
 Notice that $E_8\cong S_8^{4}$. The proof follows directly from
 Lemma \ref{lem g-non-4-multiple}, Theorems 3.6-3.9, 2.8-2.10, and
 2.16.  \qed

\vskip 0.5cm

\noindent{\bf Acknowledgments} \vskip 0.4cm

  The first author is supported by  NNSFC (Nos. 11101351 and 11171288) and NSF of the Jiangsu Higher Education Institutions (No. 11KJB110014);
  the second author is supported by  NNSFC (No. 11101351); and
  the third author is supported by  NNSFC (Nos. 11326216 and
  11301306).

\end{document}